\documentclass[reqno]{article}

\usepackage{xcolor}
\usepackage{amssymb}
\usepackage{amsfonts}
\usepackage{amsmath}
\numberwithin{equation}{section}
\usepackage{enumerate}

\newtheorem{theorem}{Theorem}[section]
\newtheorem{proposition}[theorem]{Proposition}
\newtheorem{lemma}[theorem]{Lemma}
\newtheorem{corollary}[theorem]{Corollary}
\newtheorem{definition}[theorem]{Definition}
\newtheorem{example}[theorem]{Example}

\newtheorem{remark}[theorem]{Remark}

\title {Inverse of generalized Nevanlinna function that is holomorphic at infinity}
 \author{Muhamed Borogovac, Boston, USA}
 
\begin{document}
\maketitle

\begin{abstract}

Let $\left(\mathcal{H},\left(.,.\right)\right)$ be a Hilbert space and let $\mathcal{L}\left(\mathcal{H}\right)$ be the linear space of bounded operators in $\mathcal{H}$. In this paper, we deal with $\mathcal{L}(\mathcal{H})$-valued function $Q$ that belongs to the generalized Nevanlinna class $\mathcal{N}_{\kappa} (\mathcal{H})$, where $\kappa$ is a non-negative integer. It is the class of functions meromorphic on   $C \backslash R$, such that $Q(z)^{*}=Q(\bar{z})$ and the kernel $\mathcal{N}_{Q}\left( z,w \right):=\frac{Q\left( z \right)-{Q\left( w \right)}^{\ast }}{z-\bar{w}}$ has $\kappa$ negative squares. A focus is on the functions $Q \in \mathcal{N}_{\kappa} (\mathcal{H})$ which are holomorphic at $ \infty$. A new operator representation of the inverse function $\hat{Q}\left( z \right):=-{Q\left( z \right)}^{-1}$ is obtained under the condition that the derivative at infinity $Q^{'}\left( \infty\right):=\lim\limits_{z\to \infty}{zQ(z)}$ is boundedly invertible operator. It turns out that $\hat{Q}$ is the sum $\hat{Q}=\hat{Q}_{1}+\hat{Q}_{2},\, \, \hat{Q}_{i}\in \mathcal{N}_{\kappa_{i}}\left( \mathcal{H} \right)$ that satisfies $\kappa_{1}+\kappa_{2}=\kappa $. That decomposition enables us to study properties of both functions, $Q$ and $\hat{Q}$, by studying the simple components $\hat{Q}_{1}$ and $\hat{Q}_{2}$. 
\end{abstract}
\textbf{Keywords}: Generalized Nevanlinna function, Pontryagin space, Operator representation, Generalized pole.
\\
\\
\textbf{MSC}: 47B50, 47A56, 30E99.

\section{Preliminaries and introduction }\label{s2}
\textbf{1.1 }Generalized Nevanlinna class, denoted by $\mathcal{N}_{\kappa }\left( \mathcal{H} \right)$, is extensively studied class of complex functions. For example, Hermitian matrix polynomials and their inverse functions belong to $\mathcal{N}_{\kappa }\left( \mathcal{H} \right)$. For more examples one can see, for example \cite{Lu3}. 

As usually, $N$, $R$, $C$ and $C^{+}$ denote sets of positive integers, real numbers, complex numbers, and complex numbers from the upper half plane, respectively. 

\begin{definition}\label{definition2} An operator valued complex function $Q:\mathcal{D}\left( Q \right)\to \mathcal{L}(\mathcal{H})$ belongs to the class of generalized Nevanlinna functions $\mathcal{N}_{\kappa }\left( \mathcal{H} \right)$ if it satisfies the following requirements: 

\begin{itemize}
\item $Q$ is meromorphic in $C\thinspace \backslash \thinspace R$,
\item ${Q\left( z \right)}^{\ast }=Q\left( \bar{z} \right),\thinspace z\in \mathcal{D}\left( Q \right),$ 
\item Nevanlinna kernel 
\end{itemize}
\[
\mathcal{N}_{Q}\left( z,w \right):=\frac{Q\left( z \right)-{Q\left( w \right)}^{\ast }}{z-\bar{w}},\thinspace \mathcal{N}_{Q}\left( z,\bar{z} \right):=Q^{'}\left( z \right);\thinspace z,\thinspace w\in \mathcal{D}(Q)\cap C^{+},
\]
has $\kappa $ negative squares, i.e. for arbitrary $n\in N,\thinspace  
z_{1},\thinspace \mathellipsis ,\thinspace z_{n}\in \mathcal{D}(Q)\cap 
C^{+}$ and $h_{1},\thinspace \mathellipsis ,h_{n}\in \mathcal{H}$ the Hermitian matrix 
$\left( \mathcal{N}_{Q}\left( z_{i},z_{j} \right)h_{i},h_{j} \right)_{i,j=1}^{n}$ has 
at most $\kappa $ negative eigenvalues, and for at least one choice of 
$n;\thinspace z_{1},\thinspace \mathellipsis ,\thinspace \thinspace z_{n}$, 
and $h_{1},\thinspace \mathellipsis ,h_{n}$ it has exactly $\kappa $ 
negative eigenvalues. 
\end{definition}

A generalized Nevanlinna function $Q \in \mathcal{N}_{\kappa }\left( \mathcal{H} \right)$  is called \textit{regular} if there exists at least one point $w_{0}\in \mathcal{D}(Q) \cap C^{+}$ such that the operator $Q(w_{0})^{-1}$ is boundedly invertible. 

Let $\kappa \in N\cup \left\{ 0 \right\}$ and let $\left( \mathcal{K},\, \left[ .,.\right] \right)$ denote a \textit{Krein space}. That is a complex vector space on which a scalar product, i.e. a Hermitian sesquilinear form $\left[ .,. \right]$, is defined such that the following decomposition of $\mathcal{K}$ exists
\[
\mathcal{K}=\mathcal{K}_{+} \dot{+} \mathcal{K}_{-},
\]
where $\left(\mathcal{K}_{+},\left[ .,. \right] \right)$ and $\left(\mathcal{K}_{-},-\left[ .,. \right] \right)$ are Hilbert spaces which are mutually orthogonal with respect to the form $\left[ .,. \right]$. Every Krein space $\left( \mathcal{K},\, \left[ .,.\right] \right)$ is \textit{associated} with a Hilbert space $\left( \mathcal{K},\, \left( .,.\right) \right)$, which is defined as a direct and orthogonal sum of the Hilbert spaces $\left(\mathcal{K}_{+},\left[ .,. \right] \right)$ and $\left(\mathcal{K}_{-},-\left[ .,. \right] \right)$. Topology in a Krein space $\mathcal{K}$ is introduced by means of the associated Hilbert space $\left( \mathcal{K},\, \left( .,.\right) \right)$. For properties of Krein spaces one can see e.g. \cite[Chapter V]{Bog}. 

If the scalar product $\left[ .,. \right]$ has $\kappa \, (<\infty )$ negative squares, then we call it a \textit{Pontryagin space of index} $\kappa $. The definition of a Pontryagin space and other related concepts can be found e.g. in \cite{IKL}.
\\

1.2 The following definitions of a linear relation and basic concepts related to it can be found in \cite{A,S}. In the sequel, $\mathcal{H}$, $\mathcal{K}$, $\mathcal{M}$ are inner product spaces. 

A \textit{linear relation} from $\mathcal{H}$ into $\mathcal{K}$ is a linear manifold $T$ of the product space $\mathcal{H}\times \mathcal{K}$. If $\mathcal{H}=\mathcal{K}$, $T$ is said to be a \textit{linear relation in} $\mathcal{K}$. We will use the following concepts and notations for linear relations, $T$ and $S$ from $\mathcal{H}$ into $\mathcal{K}$ and a linear relation $R$ from $\mathcal{K}$ into $\mathcal{M}$. 
\[
D\left( T \right):=\left\{ f\in \mathcal{H}\vert \left\{ f,g \right\}\in T\, for\, some\, g\in \mathcal{K} \right\},
\]
\[
R\left( T \right):=\left\{ g\in \mathcal{K}\vert \left\{ f,g \right\}\in T\, for\, some\, f\in \mathcal{H} \right\},
\]
\[
\ker T:=\left\{ f\in \mathcal{H}\vert \left\{ f,0 \right\}\in T \right\},
\]
\[
T(0):=\left\{ g\in \mathcal{K}\vert \left\{ 0,g \right\}\in T \right\},
\]
\[
T\left( f \right):=\left\{ g\in \mathcal{K}\vert \left\{ f,g \right\}\in T \right\},\, \, (f\in D\left( T \right)), 
\]
\[
T^{-1}:=\left\{ \left\{ g,f \right\}\in \mathcal{K} \times \mathcal{H}\vert \left\{ f,g \right\}\in T \right\},
\]
\[
zT:=\left\{ \left\{ f,zg \right\}\in \mathcal{H} \times \mathcal{K}\vert \left\{ f,g \right\}\in T \right\},\, \, (z\in C),
\]
\[
S+T:=\left\{ \{f,g+k\}\vert \left\{ f,g \right\}\in S,\left\{ f,k \right\}\in T \right\},
\]
\[
RT:=\left\{ \{f,k\}\in \mathcal{H} \times \mathcal{M}\vert \left\{ f,g \right\}\in T,\left\{ g,k \right\}\in R\, for\, some\, g\in \mathcal{K} \right\},
\]
\[
T^{+}:=\left\{ \{k,h\}\in \mathcal{K} \times \mathcal{H} \vert \left[ k,g \right]=\left( h,f \right)\, for\, all\, \left\{ f,g \right\}\in T \right\},
\]
\[
T_{\infty }:=\left\{ \left\{ 0,g \right\}\in T \right\}.
\]
A linear relation is \textit{closed} if it is a closed subset in the product space $\mathcal{H} \times \mathcal{K}$. If $T(0)=\lbrace 0\rbrace$, we say that $T$ is an \textit{operator}, or \textit{single-valued} linear relation.

Note, in definition of the adjoint linear relation $T^{+}$, we use the following notation for inner product spaces $\left( \mathcal{H}, (.,.) \right)$ and $\left( \mathcal{K}, [.,.] \right)$. 

Let $A$ be a linear relation in $\mathcal{K}$. We say that $A$ is \textit{symmetric} (\textit{self-adjoint}) if it holds $A\subseteq A^{+}$ ($A=A^{+})$. Every point $\alpha \in C$ for which $\left\{ f,\alpha f \right\}\in A$, with some $f\ne 0$, is called a \textit{finite eigenvalue}. The corresponding vectors are \textit{eigenvectors} belonging to the eigenvalue $\alpha $. A set that consists of all points $z\in C$ for which the relation $\left( A-zI \right)^{-1}$ is an operator defined on the entire $\mathcal{K}$, is called the \textit{resolvent} set $\rho (A)$.

It is convenient to deal with the following representation of generalized Nevanlinna functions.

\begin{theorem}\label{theorem2} A function $Q:\mathcal{D}(Q)\to \mathcal{L}(\mathcal{H})$ is a generalized Nevanlinna function of some index $\kappa $, denoted by $Q\in \mathcal{N}_{\kappa }(\mathcal{H})$, if and only if it has a representation of the form
\begin{equation}
\label{eq2}
Q\left( z \right)={Q(z_{0})}^{\ast }+(z-\bar{z_{0}})\Gamma 
_{z_{0}}^{+}\left( I+\left( z-z_{0} \right)\left( A-z \right)^{-1} 
\right)\Gamma_{z_{0}},\thinspace z\in \mathcal{D}\left( Q \right),
\end{equation}
where, $A$ is a self-adjoint linear relation in some Pontryagin space $(\mathcal{K}, [.,.])$ of index $\tilde{\kappa }\ge \kappa ; \thinspace \Gamma_{z_{0}}:\mathcal{H}\to \mathcal{K}$ is a bounded operator; $z_{0}\in \rho \left( A 
\right)\cap \mathbf{C}^{\mathbf{+}}$ is a fixed point of reference. (Then, obviously $\rho(A)\subseteq \mathcal{D}(Q)$.)  This representation can be chosen to be minimal, that is
\begin{equation}
\label{eq3}
\mathcal{K}=c.l.s.\left\{ \Gamma_{z}h:z\in \rho \left( A \right),h\in H \right\},
\end{equation}
where
\begin{equation}
\label{eq4}
\Gamma_{z}=\left( I+\left( z-z_{0} \right)\left( A-z \right)^{-1} 
\right)\Gamma_{z_{0}}.
\end{equation}

If realization (\ref{eq2}) is minimal, then $Q\in \mathcal{N}_{\kappa }(\mathcal{H})$ if and only if the negative index of the Pontryagin space $\tilde{\kappa }$ equals $\kappa $. In the case of minimal representation $\rho(A)= \mathcal{D}(Q)$ and the triple $(\mathcal{K},\thinspace A,\thinspace 
\Gamma_{z_{0}})$ is uniquely determined (up to isomorphism).
\end{theorem}
Such operator representations were developed by M. G. Krein and H. Langer, see e.g. \cite{KL1, KL2} and later converted to representations in terms of linear relations (multivalued operators), see e.g. \cite{DLS, HSW}. 

In this note, a point $\alpha \in C$ is called a \textit{finite generalized pole} of $Q$ if it is an eigenvalue of the representing relation $A$ in the minimal representation (\ref{eq2}). It means that it may be isolated singularity, i.e. an ordinary pole, as well as an embedded singularity of $Q$. The latter may be the case only if $\alpha \in R$.
\\

1.3 In this paper, we focus on the class of functions $Q\in \mathcal{N}_{\kappa }(\mathcal{H})$ that are holomorphic at $\infty$, i.e. there exists 
\begin{equation}
\label{eq6}
Q^{'}\left( \infty\right):=\lim\limits_{z\to \infty}{zQ(z)}.
\end{equation}
That is equivalent to
\begin{equation}
\label{eq8}
Q\left( z \right)=\Gamma^{+}\left( A-z \right)^{-1}\Gamma,
\end{equation}
where $A$ is a bounded self-adjoint operator in some Pontrjagin space $\mathcal{K}$, and $\Gamma:\mathcal{H}\to \mathcal{K}$ is a bounded operator, see Lemma \ref{lemma3} below. We also assume that drivative $Q^{'}\left( \infty\right)$ is boundedly invertible. In this study, $\lim\limits_{z\to \infty}{zQ(z)}$ refers to convergence in the Banach space of bounded operators $\mathcal{L}(\mathcal{H})$. By $z \to \infty$ we denote the limit if $Q$ is holomorphic at $\infty$, and  by $z \hat\to \infty$  we denote the non-tangential limit, which we use if singularities of $Q$ exist (on the real axis) in every neighborhood of $\infty$, see \cite{KL2}. The same convention applies to limits toward finite points in complex plane.   

The following well known decomposition easily follows from \cite[Proposition 3.3]{DL} for matrix functions. See \cite[Section 5.1]{Lu4} for operator valued functions.

\begin{lemma}\label{lemma2}If $Q\in \mathcal{N}_{\kappa }(\mathcal{H})$ and $\alpha \, $ is a finite generalized pole of $Q$, then it holds
\begin{equation}
\label{eq10}
{Q\left( z \right)=Q}_{\alpha }(z){+H}_{\alpha }(z),
\end{equation}
where $Q_{\alpha }\in \mathcal{N}_{\kappa_{1}}(\mathcal{H})$ is holomorphic at $\infty , H_{\alpha }\in \mathcal{N}_{\kappa 
_{2}}(\mathcal{H})$ is holomorphic at $\alpha , \kappa_{1}+\kappa_{2}=\kappa $. Then $Q_{\alpha }$ admits representation
\[
Q_{\alpha }(z)=\Gamma_{\alpha }^{+}\left( A_{\alpha }\, -z 
\right)^{-1}\Gamma_{\alpha },
\]
with a bounded operator $A_{\alpha } $. Operator $A_{\alpha }$ has the same root manifold at $\alpha $ as the representing relation $A$ of $Q$ in (\ref{eq2}).
\end{lemma}
\begin{remark}\label{remark2} The decomposition (\ref{eq10}) can be tweaked if necessary so that it holds
\[
Q(z)=\tilde{Q}(z)+\tilde{H}(z),
\]
where $\tilde{Q}(z)=\Gamma^{+}\left( \tilde{A}-z \right)^{-1}\Gamma\in 
\mathcal{N}_{\kappa_{1}}(\mathcal{H})$, self-adjoint extension $\tilde{A}$ of $A_{\alpha }$ has the same root manifold at $\alpha $ as $A_{\alpha }$, and $\Gamma^{+}\Gamma$ is a boundedly invertible operator. Then the equality  $\kappa=\kappa_{1}+\kappa_{2}$ does not have to be preserved because the number of negative squares of $\tilde{H}(z)$ may be greater than the number of negative squares of $H_{\alpha }(z).$
\end{remark} 
Indeed, if $\Gamma_{\alpha }^{+}\Gamma_{\alpha }$ is not already boundedly invertible operator in decomposition (\ref{eq10}) of $Q$, then we can add the term $\frac{B}{\beta -z}$ to $Q_{\alpha }(z)$ , where $B$ is a positive operator, $\Gamma_{\alpha }^{+}\Gamma_{\alpha }+B$ is boundedly invertible operator and $\beta \in R\setminus \lbrace \alpha \rbrace$. Also we will subtract the same term from $H_{\alpha }(z)$. Functions $\tilde{Q}\left( z \right):=Q_{\alpha }(z)+\frac{B}{\beta -z}$ and $\tilde{H}\left( z \right):=H_{\alpha }\left( z \right)-\frac{B}{\beta -z}$, will have claimed properties. $\square $
\\

1.4 The following is the summary of the main results of the paper.

In Proposition \ref{proposition6} we prove that function $Q$, which is holomorphic at $\infty $ and has invertible operator $Q^{'}(\infty)$, has $ \ker Q:= \bigcap\limits_{z\in D\left( Q \right)} \ker {Q\left( z \right)}=\lbrace  0 \rbrace$ .

The task of finding representation of $\hat{Q}(z):=-{Q\left( z \right)}^{-1}$ in terms of representing relation $A$ of $Q$ has been studied in several papers, see e.g. \cite{LaLu,Lu1}. In Theorem \ref{theorem4}, we give an operator representation of $\hat{Q}$, when function $Q$ is holomorphic at infinity and $Q^{'}\left( \infty\right)$ is boundedly invertible operator. According to Remark \ref{remark2}, those assumptions do not restrict generality in research of local properties of the function $Q \in \mathcal{N}_{\kappa }(\mathcal{H})$. 

Theorem \ref{theorem4} enables us to prove many properties of $\hat{Q}$ and $Q$. For example, in Theorem \ref{theorem6} we prove that function $Q$ which is holomorphic at $\infty$ and has $Q^{'}(\infty)$ boundedly invertible, is a regular function. In Proposition \ref{proposition12} we prove that for such $Q$ the inverse function $\hat{Q}$ must have a pole at $\infty $. In Theorem \ref{theorem10} we prove that $\hat{Q}(z):=-{Q\left( z \right)}^{-1}$ is the sum $\hat{Q}=\hat{Q}_{1}+\hat{Q}_{2},\, \, \hat{Q}_{i}\in \mathcal{N}_{\kappa_{i}}\left( \mathcal{H} \right)$, where both functions $\hat{Q}_{i}$ are represented in terms of the representing operator $A$ of $Q$, and it holds $\kappa_{1}+\kappa_{2}=\kappa $. One of the functions, say $\hat{Q}_{1}$, is a polynomial of degree one, and $\hat{Q}_{2}$ has representation of the form (\ref{eq8}). Therefore, we can call functions $\hat{Q}_{1}$ and $\hat{Q}_{2}$, \textit{polynomial}, and \textit{resolvent part of} $\hat{Q}$, respectively. Negative index $\kappa_{1}$ of $\hat{Q}_{1}$ is equal to the number of negative eigenvalues of the self-adjoint operator $\Gamma^{+}\Gamma=-\lim \limits_{z\to \infty}{zQ\left( z \right)}$. The set of zeros of $Q$ coincides with the set of poles of $\hat{Q}_{2}$. 

In Example \ref{example2}, we show how the above results can be applied to find representing operators $A$ and $\Gamma $ of $Q$ in some cases. In Example \ref{example4}, we show how to implement formulae given in Theorem \ref{theorem10} to a concrete function $Q$, in order to obtain a decomposition  $\hat Q =\hat Q_{1}+ \hat Q_{2}$  with nice properties described in that theorem.  

\section{Representation $Q\left( z \right)=S+\Gamma^{+}\left( A-z
\right)^{-1}\Gamma$}\label{s4}

2.1 We will frequently need the following proposition in this paper. 
 
\begin{proposition}\label{proposition2} 
\begin{enumerate}[(i)]
\item Let function $Q\in \mathcal{N}_{\kappa }(\mathcal{H}) $ be represented by a self-adjoint linear relation $A$ in representation (\ref{eq2}), which is not necessarily minimall. If for any point $z_{0}\in \rho\left( A \right) $ it holds 
\begin{equation}
\label{eq22}
R\left( \Gamma_{z_{0}} \right)\subseteq D\left( A \right),
\end{equation}
then the same inclusion holds for every $z\in \rho (A) $. We can define linear 
\textbf{relation} 
\begin{equation}
\label{eq24}
\Gamma:=\left( A-z \right)\Gamma_{z}, z\in \rho\left( A \right),
\end{equation}
that satisfies $D(\Gamma)=\mathcal{H}$, $\Gamma(0)=A(0)$. Then function $Q$ has representation of the form
\begin{equation}
\label{eq26}
Q\left( z \right)=S+\Gamma^{+}\left( A-z \right)^{-1}\Gamma\in 
\mathcal{N}_{\kappa }\left( \mathcal{H} \right), \, \, {S=S}^{\ast }\in \mathcal{L}\left( \mathcal{H} \right).
\end{equation}
\item Conversely, if $A$ in representation (\ref{eq26}) of $Q$ is a self-adjoint 
\textbf{linear relation} in Pontryagin space $\mathcal{K}$, and $\Gamma\subseteq \mathcal{H}\times \mathcal{K}$, $D(\Gamma)=\mathcal{H}$, is a linear relation that satisfies $A\left( 0 \right)=\Gamma(0)$, then for any point $z_{0}\in \rho(A)$ and operator 
\begin{equation}
\label{eq28}
\Gamma_{z_{0}}:=\left( A-z_{0} \right)^{-1}\Gamma,
\end{equation}
function $Q$ satisfies (\ref{eq2}). 

\item It holds 
\begin{equation}
\label{eq30}
\Gamma_{z}:=\left( I+\left( z-z_{0} \right)\left( A-z 
\right)^{-1} \right)\Gamma_{z_{0}}=\left( A-z \right)^{-1}\Gamma, \forall z \in \rho(A).\
\end{equation}
Representation (\ref{eq2}) is minimal if and only if representation (\ref{eq26}) is 
minimal.
\end{enumerate}
\end{proposition}
Note, case $S=0$ is not excluded in Proposition \ref{proposition2}. 

\textbf{Proof.} (i) For function $Q$ given by (\ref{eq2}), it holds 
\[
\Gamma_{z}=\left( I+\left( z-w \right)\left( A-z \right)^{-1} \right)\Gamma_{w}, \forall z, w \in \rho(A), 
\]
see the proof in \cite{DLS}, which obviously can be repeated when $Q \in N_{\kappa}( \mathcal H)$. If we substitute $w$ by $z_{0}$ in the above equation, then from assumption (\ref{eq22}) it follows 
\[
R\left( \Gamma_{z} \right)\subseteq D\left( A \right), \forall z \in \rho(A).
\]
In the following few steps we use properties of linear relations listed in \cite[Theorem 1.2]{A}. Note, $\Gamma_{z}$ are single-valued linear relations defined on the entire $\mathcal{H}$ which simplifies verification of the following steps. Therefore
\[
\left( A-z \right)\left( \Gamma_{z_{0}}+\left( z-z_{0} \right)\left( A-z 
\right)^{-1}\Gamma_{z_{0}} \right)=\left( A-z \right)\Gamma_{z}.
\]
According to $\left( A-z \right)\left( A-z \right)^{-1}\supseteq I$ it holds
\[
\left( A-z \right)\Gamma_{z_{0}}+\left( z-z_{0} \right)\Gamma 
_{z_{0}}\subseteq \left( A-z \right)\Gamma_{z}
\]
\[
\Rightarrow \left( A-z_{0} \right)\Gamma_{z_{0}}\subseteq \left( A-z 
\right)\Gamma_{z}, \forall z\in \rho\left( A \right).
\]
By the same token, the converse inclusion $\left( A-z \right)\Gamma 
_{z}\subseteq \left( A-z_{0} \right)\Gamma_{z_{0}}, \forall z \in \rho\left( A \right)$ holds. Therefore, 
\[
\left( A-z \right)\Gamma_{z}=\left( A-z_{0} \right)\Gamma_{z_{0}},\, 
\forall z\in \rho\left( A \right),
\]
and we can define linear relation $\Gamma$ by (\ref{eq24}). According to (\ref{eq24}) it holds 
$\Gamma(0)=A(0)$, and therefore $\left( A-z \right)^{-1}\Gamma$ is also an operator, $\forall z\in \rho\left( A \right)$.
\\

Thus, $\Gamma$ is an invariant of $Q$, i.e. $\Gamma$ is a characteristic of the function $Q$ (independent of $z\in \rho\left( A \right)$). That makes relation $\Gamma$ and representation (\ref{eq8}) particularly interesting. 



Let us now show that linear relation $\Gamma^{+}$ is an operator. If we assume the contrary, then it holds 
\[
\left\{ 0,g \right\}\in \Gamma^{+}\Rightarrow \left[ k,0 \right]=\left( 
h,g \right),\, \forall \left\{ h,k \right\}\in \Gamma.
\]
Since $D \left( \Gamma \right)=\mathcal{H}$, it follows $g=0$. Therefore, $\Gamma^{+}$ is single-valued.

From (\ref{eq24}), for $z_{0}\in \rho \left( A \right)$, we get $\Gamma=\left( 
A-z_{0} \right)\Gamma_{z_{0}}$ and $\Gamma_{z_{0}}=\left( A-z_{0} 
\right)^{-1}\Gamma$. Then we substitute $\Gamma_{z_{0}}^{+}$ and 
$\Gamma_{z_{0}}$ into (\ref{eq2}) and easily derive
\[
Q\left( z \right)=Q\left( \bar{z_{0}} \right)+\left( z-\bar{z_{0}} 
\right)\Gamma^{+}\left( A-\bar{z_{0}} \right)^{-1}\left( A-z 
\right)^{-1}\Gamma.
\]
By means of the resolvent equation we get 
\[
Q\left( z \right)=Q\left( \bar{z_{0}} \right)-\Gamma^{+}\left( 
A-\bar{z_{0}} \right)^{-1}\Gamma+\Gamma^{+}\left( A-z 
\right)^{-1}\Gamma.
\]
By substituting here 
\[
S:=Q\left( \bar{z_{0}} \right)-\Gamma^{+}\left( A-\bar{z_{0}} 
\right)^{-1}\Gamma,
\]
we get the first equation of (\ref{eq26}). 

From the first equation of (\ref{eq26}) and from ${Q\left( z \right)}^{\ast }=Q\left( 
\bar{z} \right)$ it follows ${S=S}^{\ast }$. 

\textbf{(ii)} Conversely, assume (\ref{eq26}) holds with linear relation $A$. From 
(\ref{eq26}), for ${z=z}_{0}$, we get ${S=S}^{\ast }={Q\left( z_{0} \right)}^{\ast 
}-\Gamma^{+}\left( A-\bar{z_{0}} \right)^{-1}\Gamma$. Substituting 
$S$ into (\ref{eq26}) and applying resolvent equation we obtain
\[
Q\left( z \right)={Q\left( z_{0} \right)}^{\ast }+\left( z-\bar{z_{0}} 
\right)\Gamma^{+}{\left( A-\bar{z_{0}} \right)^{-1}\left( A-z 
\right)}^{-1}\Gamma.
\]
Now (\ref{eq28}) gives 
\[
Q\left( z \right)={Q\left( z_{0} \right)}^{\ast }+\left( z-\bar{z_{0}} 
\right)\Gamma_{z_{0}}^{+}\left( A-z \right)^{-1}\Gamma.
\]
According to resolvent equation it holds
\begin{equation}
\label{eq29}
\left( A-z \right)^{-1}=\left( I+\left( z-z_{0} \right)\left( A-z \right)^{-1} \right)\left( A-z_{0} \right)^{-1}, \forall z\in \rho(A).
\end{equation}
Therefore
\[
Q\left( z \right)={Q\left( z_{0} \right)}^{\ast }+\left( z-\bar{z_{0}} 
\right)\Gamma_{z_{0}}^{+}\left( I+\left( z-z_{0} \right)\left( A-z 
\right)^{-1} \right)\left( A-z_{0} \right)^{-1}\Gamma.
\]
Substituting here $\Gamma_{z_{0}}$ from (\ref{eq28}) gives (\ref{eq2}). 
\\

(iii) From (\ref{eq29}) and (\ref{eq28}) it follows
\[
\left( A-z \right)^{-1}\Gamma=\left( I+\left( z-z_{0} \right)\left( A-z 
\right)^{-1} \right)\Gamma_{z_{0}}=:\Gamma_{z}, \forall z\in \rho(A).\
\]
This proves (\ref{eq30}).
Minimality of a representation is defined in terms the of vectors $\Gamma_{z}h$ by (\ref{eq3}). According to (\ref{eq30}) we conclude that represention (\ref{eq26}) is minimal if and only if
\[
\mathcal{K}=c.l.s.\left\{ \left( A-z \right)^{-1}\Gamma h:z\in \rho \left( A \right),h\in 
\mathcal{H} \right\}.
\]
This proves (iii). 
$\square$
\\

Note, the first statement of the proposition is well known for matrix functions represented by operators. This was proven in \cite{KL2} for scalar, and in \cite{LaLu} for matrix valued function $Q$. In both cases one additional assumption on $Q$ was made so that $A$ was linear operator from the start. 
\\

By definition, $\infty$ is generalized pole of $Q$ if and only if $0$ is generalized pole of the function $\tilde Q\left(\zeta \right) = Q \left(\frac{-1}{\zeta}\right)$, see \cite[Remark 3.13.]{BLu}. This is equivalent to $A(0)\neq \lbrace 0\rbrace$, where $A$ is representing relation of $Q$. In that case $\infty$ is called an eigenvalue of $A$ and nonzero vectors from $A(0)$ are called \textit{eigenvectors at} $\infty$, see \cite{Lu2}.  
\\

The following statement is well known for closed linear relations in Hilbert space $\mathcal{H}$, see e.g. \cite{LT}. We will state it here in our setting, for convenience of the reader. 

\begin{lemma}\label{lemma3a} 
Let $\mathcal{H}$ and $\mathcal{K}$ be Hilbert and Krein space, respectively, and let linear relation $T\subseteq \mathcal{H}\times \mathcal{K}$ has closed $T(0)$. Then it holds:
\[
T=\tilde{T}\dot{+}T_{\infty },
\]
where $\dot{+}$ denotes direct sum of subspaces, $\tilde{T}$ is an operator with $D\left( \tilde{T} \right)=D\left( T \right)$ and $T_{\infty }:=\left\{ \left\{ 0,g \right\}\in T \right\}$.
\end{lemma}
\textbf{Proof}: Because $T(0)\subseteq \mathcal{K}$ is closed subspace of the Hilbert space $\left( \mathcal{K}, \left( .,.\right)\right) $ associated with Krein space $\left( \mathcal{K}, \left[ .,.\right]\right) $, we can uniquely and orthogonaly decompose $\left( \mathcal{K}, \left( .,.\right)\right) $ by means of $T(0)$. Thus, for every $\left\{ f,g \right\}\in T$ we have, $\left\{ f,g \right\}=\left\{ f,g_{1}(\dot{+})g_{0} \right\}$, where $(\dot{+})$ is direct and orthogonal sum in the Hilbert space $\left( \mathcal{K}, \left( .,.\right)\right)$, and $g_{0}\in T(0)$ and $g_{1}\in \mathcal{K}(-)T(0)$ are uniquely determined vectors. We define 
\[ 
\tilde{T}:=\left\{ \left\{ f,g_{1} \right\}\vert \left\{ f,g \right\}\in T \right\},
\]
and $T_{\infty}$ is as above. Then we have 
\[
T=\tilde{T}(\dot{+})T_{\infty }\subseteq \mathcal{H}\times \mathcal{K},
\]
where $(\dot{+})$ denotes direct orthogonal sum in the Hilbert space associated with $\mathcal{H}\times \mathcal{K}$. 

Because the sum $g_{1}(\dot{+})g_{0}$ does not have to be orthogonal in the Krein space $\left( \mathcal{K}, \left[ .,.\right]\right) $, we write
\[
T=\tilde{T}\dot{+}T_{\infty }.
\]
It is easy to verify that $\tilde{T}=T(-)T_{\infty }$ is single-valued.$\quad$ $\square$
\begin{corollary}\label{corollary2} If representing relation $A$ of $Q\in \mathcal{N}_{\kappa }(\mathcal{H})$  satisfies condition (\ref{eq22}), then $A$ can be replaced in (\ref{eq2}) by its operator part $\tilde{A}$. If representation (\ref{eq2}) is minimal, it will remain minimal with self-adjoint operator $\tilde{A}$. The function $Q$ does not have generalized pole at $\infty$. 
\end{corollary}
\textbf{Proof}. Because $A$ is closed linear relation, it is easy to verify that $A(0)$ is closed. According to Lemma \ref{lemma3a} it holds
\[
A=\tilde{A}\dot{+}A_{\infty }.
\]
According to Proposition \ref{proposition2} (i) there exists a linear relation 
\[
\Gamma:=\left( A-z \right)\Gamma_{z},\, z\in \rho(A),
\]
with $\Gamma\left( 0 \right)=A\left( 0 \right)$. Because $\Gamma\left( 0 \right)$ is closed, according to Lemma \ref{lemma3a} it holds 
\[
\Gamma =\tilde{\Gamma }\dot{+}\Gamma_{\infty }.
\]
Because $\Gamma\left( 0 \right)=A\left( 0 \right)=\ker \left( A-z \right)^{-1}$, it holds
\begin{equation}
\label{eq30a}
\Gamma_{z}=\left( A-z \right)^{-1}\Gamma=\left( \tilde A-z \right)^{-1}\tilde \Gamma , \forall z\in \rho(A) .
\end{equation}

Let $z_{0}\in \rho(A)\setminus R $ be the point of reference in (\ref{eq2}). Let us now prove that we can replace $\left( A-z \right)^{-1}\Gamma_{z_{0}}$ by  $\left( \tilde A-z \right)^{-1} \Gamma_{z_{0}}$ in  (\ref{eq2}). We start from (\ref{eq4}) written in the form
\[
\left( A-z \right)^{-1}\Gamma_{z_{0}}=\, \frac{\Gamma_{z}-\Gamma 
_{z_{0}}}{z-z_{0}}, \forall z \in \rho(A) .
\]
According to (\ref{eq30a}) and the resolvent equation we have 
\[
\left( A-z \right)^{-1}\Gamma_{z_{0}}=\frac{\left( \tilde{A}-z 
\right)^{-1}\tilde{\Gamma }-\left( \tilde{A}-z_{0} \right)^{-1}\tilde{\Gamma 
}}{z-z_{0}}=\left( \tilde{A}-z \right)^{-1}\left( \tilde{A}-z_{0} 
\right)^{-1}\tilde{\Gamma }=\left( \tilde{A}-z \right)^{-1}\Gamma_{z_{0}}.
\]
This proves  
\[
\left( A-z \right)^{-1}\Gamma_{z_{0}}=\left( \tilde{A}-z \right)^{-1}\Gamma_{z_{0}}.
\]

Therefore, we can substitute $\left( \tilde{A}-z \right)^{-1}\Gamma_{z_{0}}$ for $\left( A-z \right)^{-1}\Gamma_{z_{0}}$ into (\ref{eq4}) and (\ref{eq2}), and values of $\Gamma_{z}$ and $Q(z)$ will not change. Thus, 
\[
\Gamma_{z}= \left( I+\left( z-z_{0} \right)\left( \tilde{A}-z \right)^{-1} \right)\Gamma_{z_{0}}.
\] 
\[
Q\left( z \right)={Q(z_{0})}^{\ast }+(z-\bar{z_{0}})\Gamma _{z_{0}}^{+}\left( I+\left( z-z_{0} \right)\left( \tilde{A}-z \right)^{-1} \right)\Gamma_{z_{0}},\thinspace z\in \mathcal{D}\left( Q \right).
\]
According to definition of minimality (\ref{eq3}), we conclude that minimal representation (\ref{eq2}) remains minimal when $\tilde{A}$ replaces $A$. Because of the uniqueness of the minimal representation (\ref{eq2}) it must be $A=\tilde{A}$. Therefore, $\tilde A$ must be a self-adjoint operator, as the unique representing operator of a generalized Nevanlinna function. Because the function $Q$ is represented by operator $\tilde A$, we conclude that $Q$ cannot have generalized pole at $\infty$. $\square$
\\

2.2 By definition a function $Q$ has a non-tangential limit at $\infty$ if and only if the function  $\tilde Q\left(\zeta \right) = Q \left(\frac{-1}{\zeta}\right)$ has a non-tangential limit at $0$. By the same token a function $Q$ is holomorphic at $\infty$ if and only if the function  $\tilde Q\left(\zeta \right) = Q \left(\frac{-1}{\zeta}\right)$ is holomorphic at $0$. The following proposition, that corresponds to \cite[Satz 1.4]{KL2} holds. 

\begin{proposition}\label{proposition4} Let $Q\in \mathcal{N}_{\kappa }(\mathcal{H})$ satisfies non-tangential version of (\ref{eq6}):
\begin{equation}
\label{eq30b}
\exists Q^{'}\left( \infty \right):=\lim \limits_{z \hat\to \infty}{zQ(z)},
\end{equation}
where the limit denotes convergence in the Banach space of bounded operators. Then $Q^{'}(\infty)\in \mathcal{L}(\mathcal{H})$, and $Q$ has minimal representation (\ref{eq2}) with a self-adjoint operator $A$.
\end{proposition}

\textbf{Proof}. Because $\mathcal{L}(\mathcal{H})$ is a Banach space with respect to norm topology, we conclude that $Q^{'}(\infty)$, given by (\ref{eq30b}), is a bounded operator. Under assumption that limit (\ref{eq30b}) exists, it holds 
\[ 
\lim \limits_{\zeta \hat\to 0}\tilde{Q}(\zeta):=\lim \limits_{z \hat\to \infty}{Q(z)}=0.
\] 
If we define:  $\tilde{Q}(0):=\lim \limits_{\zeta \hat\to 0}{\tilde{Q}(\zeta)}=0$, then 
\[ 
\tilde{Q}^{'}(0):=\lim \limits_{\zeta \hat\to 0}\frac{\tilde{Q}(\zeta)-\tilde{Q}(0)}{\zeta}=\lim \limits_{z \hat\to \infty}{zQ(z)}=:Q^{'}(\infty).
\]
According to \cite[Defintion 3.1 (B)]{BLu}, $\zeta = 0$ is not a generalized pole of $\tilde Q$, i.e. $\infty$ is not a generalized pole of Q. Therefore, the representing relation $A$ satisfies $A(0)={0}$. Hence, $Q$ is represented by the self-adjoint operator $A$ in (\ref{eq2}).$\square$

\begin{lemma}\label{lemma3} A function $Q\in \mathcal{N}_{\kappa }(\mathcal{H})$ is holomorphic at $\infty$ if and only if $Q(z)$ has minimal representation (\ref{eq8}) 
\[
Q\left( z \right)=\Gamma^{+}\left( A-z \right)^{-1}\Gamma, \thinspace z\in \mathcal{D}(Q), 
\]
with a bounded self-adjoint operator  $A$ in a Pontryagin space $\mathcal{K}$, and bounded operator $\Gamma:\mathcal{H}\to \mathcal{K}$. In this case 
\[
Q^{'}\left( \infty \right):=\lim \limits_{z\to \infty}{zQ(z)}=-\Gamma^{+} \Gamma.
\]
\end{lemma}
\textbf{Proof.} If $Q(z)$ is holomorphic at $\infty$, then it satisfies (\ref{eq30b}). According to Proposition \ref{proposition4}, $Q$ is represented by an operator $A$. From the assumption of holomorphy at $\infty$ it follows that operator $A$ has bounded spectrum. According to \cite[Corollary 2]{L}, $A$ is bounded. Then condition (\ref{eq22}) is satisfied. According to Proposition \ref{proposition2} (i), $Q$ has minimal representation (\ref{eq26}). Then, from existence of limit (\ref{eq30b}), it follows $S=0$.

Conversely, if $A$ is bounded operator in representation (\ref{eq8}), then it has bounded spectrum, and therefore, $Q$ is holomorphic at infinity. 

To prove the last statement of the lemma, we use Neumann series of  resolvent of the bounded operator $A$. 
\[  
Q^{'}\left( \infty \right):=\lim \limits_{z\to \infty}{zQ(z)}=\lim\limits_{z\to\infty}{z\Gamma^{+}\left(\sum\limits_{i=0}^\infty {-\frac{A^{i}}{z^{i+1}}} \right)\Gamma=-\Gamma^{+}\Gamma}.
\] $\square$

The concept 
\[
\ker Q:= \bigcap\limits_{z\in D\left( Q \right)} \ker {Q\left( z \right)} 
\]
was introduced in \cite{DLS}. For matrix function $Q\in N_{\kappa }^{n \times n}$, represented by (\ref{eq2}) it 
was proven 
\[
\ker Q=\ker \Gamma_{z_{0}}\bigcap \ker {Q(z_{0})}^{\ast }. 
\]
\begin{proposition}\label{proposition6} If $Q\in N_{\kappa }(\mathcal H)$ is holomorphic at 
infinity and $Q^{'}\left( \infty \right)$ is invertible, then 
\[
\ker Q=\left\{ 0 \right\}.
\] 
\end{proposition}

\textbf{Proof}. According to Lemma \ref{lemma3} we can assume that $Q$ is minimally represented by bounded operator $A$. Recall, for $z, w \in \rho\left( A \right)=\mathcal{D}(Q)$ it holds
\[
\Gamma_{z}=\left( I+\left( z-w \right)\left( A-z \right)^{-1} \right)\Gamma_{w}. 
\] 
Obviously, 
\[
\Gamma_{w}h=0\Rightarrow \Gamma_{z}h=0,
\]
If we reverse roles of $z$ and $w$, then the converse implication holds. 
Hence, it holds
\[
\ker \Gamma_{z}=\ker \Gamma_{w}.
\] 
If $Q(z)$ is holomorphic at $\infty $, according to Lemma \ref{lemma3},  $Q$ has representation (\ref{eq8}) with bounded operator $A$. Therefore, condition (\ref{eq22}) is satisfied. According to Proposition \ref{proposition2} (iii) we have 
\[
 \Gamma_{z}=\left( A-z \right)^{-1}\Gamma , \forall z \in \mathcal{D}(Q). 
\]
Then we have:
\[
\left( \ref{eq8} \right)\Rightarrow Q\left( z \right)h=\Gamma^{+}\Gamma_{z}h, \thinspace \forall h \in \mathcal{H}, \forall z \in \mathcal{D}(Q).
\]
If we assume $ h\in \ker Q$, then according to definition of $\ker Q$ we have 
\[
h \in \ker {Q} \Leftrightarrow h \in \ker {zQ(z)}, \forall z \in \mathcal{D}(Q)
\]
\[
\Leftrightarrow 0={\lim \limits_{z \hat{\to} \infty}{zQ(z)h}} =-{\Gamma^{+}\Gamma}h=Q^{'}(\infty)h \Leftrightarrow h = 0.  
\]
This proves the statement. $ \square$
\\

We cannot here claim that $Q(z)$ is a regular function. We will prove it in the following section. 






\section{Inverse of ${\Gamma}^{+}\left(A-z\right)^{-1}\Gamma$}\label{s6}
\begin{lemma}\label{lemma8} Let bounded operators $\Gamma :\mathcal{H}\to \mathcal{K}$ and $\Gamma^{+}:\mathcal{K}\to \mathcal{H}$ be introduced as usually, see section \ref{s2}. Assume also that $\Gamma^{+}\Gamma 
$ is a boundedly invertible operator in the Hilbert space $\left( \mathcal{H},\left( .,. \right) \right)$. Then for operator 
\begin{equation}
\label{eq31}
P:=\Gamma\left( \Gamma^{+}\Gamma \right)^{-1}\Gamma^{+}
\end{equation}
the following statements hold:

\begin{enumerate}[(i)]
\item P is orthogonal projection in Pontryagin space $(\mathcal{K},\, \left[ .,. \right]).$
\item Scalar product does not degenerate on $\Gamma \left( \mathcal{H} \right)$ and therefore it does not degenerate on ${\Gamma \left( \mathcal{H} \right)\, }^{[\bot ]}=\ker \Gamma^{+}. $
\item $\ker \Gamma^{+}=\left( I-P \right)\mathcal{K}.$ 
\item Pontryagin space $\mathcal{K}$ can be decomposed as a direct orthogonal sum of Pontryagin spaces i.e. 
\end{enumerate}
\begin{equation}
\label{eq32}
\mathcal{K}=\left( I-P \right)\mathcal{K}[+]P\mathcal{K}.
\end{equation}
\end{lemma}
\textbf{Proof.} (i) Obviously $P^{2}=P$. 

According to well known properties of adjoint operators, see e.g. \cite[p.34]{IKL}, it is easy to verify $\left[ \left( \Gamma^{+}\Gamma \right)^{-1} \right]^{\ast }=\left( \Gamma^{+}\Gamma \right)^{-1}$ and then to verify $\left[ Px,y \right]=\left[ x,Py \right],\, i.e.\, P^{[\ast ]}=P$. This proves (i). 
\\

(ii) If $\Gamma h\ne 0$ and $\left[ \Gamma h,\Gamma g \right]=0$, $\forall 
g\in \mathcal{H}$, then $({\, \Gamma }^{+}\Gamma h,g)=0$, $\forall g\in \mathcal{H}$. Then we have
\\
${\, \Gamma }^{+}\Gamma h=0 \Rightarrow h=0 \Rightarrow \Gamma h=0.$ This is a contradiction that 
proves (ii).
\\

(iii) It is sufficient to prove $\ker \Gamma^{+}=\ker P$. 
\[
P:=\Gamma \left( \Gamma^{+}\Gamma \right)^{-1}\Gamma^{+}\Rightarrow 
\ker \Gamma^{+}\subseteq \ker P.
\]
Conversely, because $\Gamma^{+}\Gamma $ is boundedly invertible $R\left( \Gamma 
^{+} \right)=\mathcal{H}$. Then
\[
y\in \ker P\Rightarrow 0=\left[ \Gamma \left( \Gamma^{+}\Gamma 
\right)^{-1}\Gamma^{+}y,x \right]=\left( \left( \Gamma^{+}\Gamma 
\right)^{-1}\Gamma^{+}y,\Gamma^{+}x \right),\forall \Gamma^{+}x\in \mathcal{H}.
\]
\[
R\left( \Gamma^{+} \right)=\mathcal{H}\Rightarrow \left( \Gamma^{+}\Gamma 
\right)^{-1}\Gamma^{+}y=0\Rightarrow \Gamma^{+}y=0\Rightarrow y\in \ker 
\Gamma^{+}.
\]
(iv) This statement follows directly from (iii) and (ii). \quad $\square$
\\

Assume now that function $Q$ is given by (\ref{eq8}) and that projection $P$ is given by (\ref{eq31}). We define  
\[
\tilde{A}:=\left( I-P \right){A}_{\mid( I-P )\mathcal{K}}\, .
\]
Then
\[
{(\tilde{A}-zI_{\mid( I-P )\mathcal{K}})}^{-1}:\left( I-P \right)\mathcal{K} \rightarrow \left( I-P \right)\mathcal{K}.
\]
Note that it is customary to omit the identity mapping in resolvents. Therefore, we will frequently write $\left( \tilde{A}-z \right)^{-1}$ rather than $\left( \tilde{A}-zI_{\mid( I-P )\mathcal{K}} 
\right)^{-1}$. It holds
\[
\left( I-P \right)\left( \tilde{A}-z \right)^{-1}\left( I-P \right)=\left( 
{\begin{array}{*{20}c}
\left( \tilde{A}-zI_{\mid( I-P )\mathcal{K}} \right)^{-1} & 0\\
0 & 0\\
\end{array} } \right).
\]
In the sequel, we will use notation from the left hand side of this equation 
because it makes the following proofs easier to write. 
\\

\begin{theorem}\label{theorem4}Assume that function $Q\in \mathcal{N}_{\kappa}(\mathcal{H})$ is holomorphic at $\infty $, and that
\[
Q^{'}\left( \infty\right):=\lim\limits_{z\to \infty}{zQ(z)}
\]
is boundedly invertible. Then there exists the inverse function 
\[
\hat{Q}\left( z \right):=-{Q(z)}^{-1},\, 
\]
and $\hat{Q}\left( z \right)$ has the following representation on $\mathcal{D}(Q)\cap \mathcal{D}(\hat Q)$
\begin{equation}
\label{eq36}
\hat{Q}\left( z \right)=\left( \Gamma^{+}\Gamma \right)^{-1}\Gamma^{+}\left\{ A(I-P)\left( \tilde{A}-z \right)^{-1}(I-P)A-\left( A-zI \right) \right\}\Gamma\left( \Gamma^{+}\Gamma \right)^{-1},
\end{equation}
where operator $\Gamma$ was defined by (\ref{eq24}) and projection $P$ was defined by equation (\ref{eq31}). 
\end{theorem}

\textbf{Proof.} According to Lemma \ref{lemma3}, function $Q$ has minimal representation 
(\ref{eq8}) with bounded operator $A$. For projection $P$ defined in Lemma \ref{lemma8}, we have the following decomposition with respect to (\ref{eq32}) 
\[
A-zI=\left( {\begin{array}{*{20}c}
\left( I-P \right)(A-zI)(I-P) & \left( I-P \right)AP\\
PA(I-P) & P\left( A-zI \right)P\\
\end{array} } \right).
\]
Let us denote 
\[
\left( {\begin{array}{*{20}c}
X & Y\\
Z & W\\
\end{array} } \right):=\left( A-z \right)^{-1}\, .
\]
By solving operator equations derived from the identity 
\[
\left( {\begin{array}{*{20}c}
X & Y\\
Z & W\\
\end{array} } \right)\left( {\begin{array}{*{20}c}
\tilde{A}-z(I-P) & \left( I-P \right)AP\\
PA(I-P) & P\left( A-zI \right)P\\
\end{array} } \right)=\left( {\begin{array}{*{20}c}
I-P & 0\\
0 & P\\
\end{array} } \right)
\]
we get 
\[
W=\left\{ P\left( A-zI \right)P-PA(I-P)\left( \tilde{A}-z 
\right)^{-1}(I-P)AP \right\}^{-1}.
\]
It is easy to verify the following equalities:
\[
\Gamma^{+}P=\Gamma^{+},\, \, P\Gamma=\Gamma,\, \, \Gamma^{+}\left( I-P \right)=0,\, \, \left( I-P \right)\Gamma=0.
\]
It follows 
\[
Q\left( z \right)=\Gamma^{+}\left( {\begin{array}{*{20}c}
X & Y\\
Z & W\\
\end{array} } \right)\Gamma=\, \left( \Gamma^{+}\left( I-P 
\right),\Gamma^{+}P\, \right)\left( {\begin{array}{*{20}c}
X & Y\\
Z & W\\
\end{array} } \right)\left( {\begin{array}{*{20}c}
\left( I-P \right)\Gamma\\
P\Gamma\\
\end{array} } \right)
\]
\[
\Rightarrow Q\left( z \right)=\left( 0,\Gamma^{+}\, \right)\left( 
{\begin{array}{*{20}c}
X & Y\\
Z & W\\
\end{array} } \right)\left( {\begin{array}{*{20}c}
0\\
\Gamma\\
\end{array} } \right)=\Gamma^{+}\left( {\begin{array}{*{20}c}
0 & 0\\
0 & W\\
\end{array} } \right)\Gamma.
\]
Therefore, we do not need to find operators $X$, $Y$, $Z$. By substituting $W$ here, we get
\begin{equation}
\label{eq40}
Q\left( z \right)=\Gamma^{+}\left\{ P\left( A-zI \right)P-PA(I-P)\left( 
\tilde{A}-z \right)^{-1}(I-P)AP \right\}^{-1}\Gamma.
\end{equation}
By substituting expressions (\ref{eq40}) and (\ref{eq36}) for $Q$ and $\hat{Q}$, 
respectively, into the following product, we verify 
\[
Q\left( z \right)\hat{Q}\left( z \right)=
\]
\[
=\Gamma^{+}\left\{ P\left( A-zI \right)P-PA(I-P)\left( \tilde{A}-z 
\right)^{-1}(I-P)AP \right\}^{-1}\Gamma\times
\]
\[
\left( \Gamma^{+}\Gamma \right)^{-1}\Gamma^{+}\left\{A(I-P)\left( \tilde{A}-z \right)^{-1}(I-P)A-\left( A-zI \right) 
\right\}\Gamma\left( \Gamma^{+}\Gamma \right)^{-1}=
\]
\[
=\Gamma^{+}\left\{ P\left( A-zI \right)P-PA(I-P)\left( \tilde{A}-z 
\right)^{-1}(I-P)AP \right\}^{-1}\times
\]
\[
\left\{ PA(I-P)\left( \tilde{A}-z \right)^{-1}(I-P)AP-P\left( A-zI \right)P 
\right\}\Gamma\left( \Gamma^{+}\Gamma \right)^{-1}
\]
\[
=\Gamma^{+}\left( -P \right)\Gamma\left( \Gamma^{+}\Gamma \right)^{-1}=-I.
\]$\square$

The remaining statements of this paper are consequences of Theorem \ref{theorem4}. 

\begin{theorem}\label{theorem6} Let $Q\in \mathcal{N}_{\kappa }(\mathcal{H})$. 
\begin{enumerate}[(i)]
\item $Q$ is holomorphic at $\infty $ and $Q^{'}\left( \infty\right)$ is boundedly invertible if and only if 
\begin{equation}
\label{eq42}
\hat{Q}\left( z \right)=\tilde{\Gamma }^{+}\left( \tilde{A}-z \right)^{-1}\tilde{\Gamma }+\hat{S}+\hat{G}z, \forall z \in \mathcal{D}(Q)\cap \mathcal{D}(\hat Q)
\end{equation}
where $\tilde{A}$ is a self-adjoint bounded operator in the Pontryagin space $(I-P)\mathcal{K}$, $\hat{S}$ and $\hat{G}$ are self-adjoint bounded operators in the Hilbert space $\mathcal{H}$, and $\hat{G}$ is boundedly invertible. 
\item In that case function $Q\in \mathcal{N}_{\kappa }(\mathcal{H})$ is regular.
\end{enumerate}
\end{theorem}

\textbf{Proof} (i) $(\Rightarrow)$ The assumptions are the same as in Theorem \ref{theorem4}. Therefore, representation (\ref{eq36}) holds. If we substitute 
\begin{equation}
\label{eq42a}
\hat{S}=-\left( \Gamma^{+}\Gamma \right)^{-1}\Gamma^{+}A\Gamma\left( \Gamma^{+}\Gamma \right)^{-1},\, 
\hat{G}=\left( \Gamma^{+}\Gamma \right)^{-1}
\end{equation}
\begin{equation}
\label{eq42b}
\tilde{\Gamma }:=\left( I-P \right)A\Gamma \left( \Gamma^{+}\Gamma \right)^{-1},
\end{equation}
into representation (\ref{eq36}) we get representation (\ref{eq42}). Operator $\tilde A$ is bounded because it is a restriction of the bounded operator $A$. The statements about $\hat S$ and $\hat G$ are easy verification.
\\

$(\Leftarrow)$ Now we assume that (\ref{eq42}) holds. Obviously:
\[
\lim\limits_{z \to \infty}\frac{\hat{Q}\left( z \right)}{z}=\lim\limits_{z \to \infty}\left(-zQ(z) \right)^{-1}.
\]
On the other hand, because $\tilde{A}$ is bounded we can apply Neumann 
series of the resolvent $\left( \tilde{A}-z \right)^{-1}$. We have 
\[
\lim\limits_{z \to \infty}\frac{\hat{Q}\left( z \right)}{z}=\lim\limits_{z \to \infty}\left( \frac{\tilde{\Gamma }^{+}\left(\tilde{A}-z \right)^{-1}\tilde{\Gamma }+\hat{S}}{z}+\hat{G} \right)=
\]
\[
\lim\limits_{z \to \infty} \left( \tilde{\Gamma }^{+}\sum\limits_{i=0}^\infty 
{-\frac{\tilde{A}^{i}}{z^{i+2}}} \tilde{\Gamma }+\frac{\hat{S}}{z}\right)+\hat{G}=\hat{G}.
\]
Therefore, 
\[
\lim\limits_{z \to \infty}\left( -zQ(z) \right)^{-1}=\hat{G}.
\]
Because $\hat{G}$ is bounded, $\lim\limits_{z\to \infty}{zQ\left( z \right)}$ is boundedly invertible.

(ii) This statement holds because, according to (\ref{eq42}), operator $\hat{Q}(z)$ is obviously bounded for every $z \in \mathcal{D}(Q)\cap \mathcal{D}(\hat Q)$ . $ \square$ 
\\

It is usually very difficult to find representing operator for a given function $Q\in \mathcal{N}_{\kappa }(\mathcal{H})$. The construction used in cited papers is abstract and not applicable in concrete situations. Theorem \ref{theorem4} gives us a new simple relationships between representing operators $A$, $\Gamma$ and $\Gamma^{+}$. That might help us to find those operators in some cases, like e.g. in the following case.

\begin{example}\label{example2} Given function
\[
Q\left( z \right)=-\left[ {\begin{array}{*{20}c}
0 & z^{-1}\\
z^{-1} & z^{-2}\\
\end{array} } \right].
\]
\end{example}
It is easy to verify that function $Q(z)$ is holomorphic at infinity, and that it holds 
\[
Q^{'}(\infty):=\lim\limits_{z \to \infty}{zQ\left( z \right)}=-\left[ 
{\begin{array}{*{20}c}
0 & 1\\
1 & 0\\
\end{array} } \right].
\]
According to Lemma \ref{lemma3},  $Q(z)$ admits minimal representation (\ref{eq8}). Hence, 
\[
Q\left( z \right)=\Gamma^{+}\left( A-zI \right)^{-1}\Gamma \wedge -\left[ 
{\begin{array}{*{20}c}
0 & 1\\
1 & 0\\
\end{array} } \right]=-\Gamma^{+}\Gamma.
\]
In addition, 
\[
{Q\left( z \right)}^{-1}=\left[ {\begin{array}{*{20}c}
1 & -z\\
-z & 0\\
\end{array} } \right]=:L\left( z \right).
\]
i.e. the inverse function is a polynomial. Therefore, the resolvent part of $\hat{Q}$ in representation (\ref{eq36}) must be equal to zero. It holds,
\[
\left( \Gamma^{+}\Gamma \right)^{-1}\Gamma^{+}\left( A-zI 
\right)\Gamma\left( \Gamma^{+}\Gamma \right)^{-1}=\left[ 
{\begin{array}{*{20}c}
1 & -z\\
-z & 0\\
\end{array} } \right]
\]
\[
\Rightarrow \Gamma^{+}\left( A-zI \right)\Gamma=\left[ 
{\begin{array}{*{20}c}
0 & -z\\
-z & 1\\
\end{array} } \right]
\Rightarrow  
\Gamma^{+}A\Gamma 
= \Gamma^{\ast }JA\Gamma=\left[ 
{\begin{array}{*{20}c}
0 & 0\\
0 & 1\\
\end{array} } \right].
\]
Here $J$ denotes a fundamental symmetry in $\mathcal{K}$. Because function $Q$ has a single 
pole of order two at $z=0$, the representing operator has the single 
eigenvalue of order two at $z=0$. All those information enable us to make an easy educated guess
\[
A=\left[ {\begin{array}{*{20}c}
0 & 1\\
0 & 0\\
\end{array} } \right], \Gamma=\left[ {\begin{array}{*{20}c}
1 & 0\\
0 & 1\\
\end{array} } \right],\, J=\left[ {\begin{array}{*{20}c}
0 & 1\\
1 & 0\\
\end{array} } \right]=\Gamma^{+}.
\] $\square$

We will refer to this example for a different reason in Theorem \ref{theorem10}.

\begin{proposition}\label{proposition10} Let $Q\left( z \right),\hat{Q}\left( z \right), \Gamma, {\Gamma}^{+}$ be the same as in Theorem \ref{theorem4}. Then for all $z \in \mathcal{D}(Q)\cap \mathcal{D}(\hat Q)$ it holds
\begin{equation}
\label{eq46}
\hat{Q}\left( z \right)\Gamma^{+}=\left( \Gamma^{+}\Gamma 
\right)^{-1}\Gamma^{+}\left\{ -I+A(I-P)\left( \tilde{A}-z 
\right)^{-1}(I-P) \right\}(A-zI).
\end{equation}
\end{proposition}

\textbf{Proof. }In the following derivations we will frequently use $\Gamma^{+}P=\Gamma^{+}$ and $\, P\Gamma=\Gamma$. From (\ref{eq36}) it 
follows
\[
\hat{Q}\left( z \right)\Gamma^{+}=\left( \Gamma^{+}\Gamma \right)^{-1}\Gamma^{+}\left\{A(I-P)\left( \tilde{A}-z \right)^{-1}(I-P)A-\left( A-zI \right) 
\right\}\Gamma\left( \Gamma^{+}\Gamma \right)^{-1}\Gamma^{+}
\]
\[
=\left( \Gamma^{+}\Gamma \right)^{-1}\Gamma^{+}\left\{ 
A(I-P)\left( \tilde{A}-z \right)^{-1}(I-P)(A-zI)P-\left( A-zI \right)P 
\right\}
\]
\[
=\left( \Gamma^{+}\Gamma \right)^{-1}\Gamma^{+}\lbrace 
A\left( I-P \right)\left( \tilde{A}-z \right)^{-1}\left( I-P \right)\left(
A-zI \right)\left( P-I \right)+
\]
\[
+A(I-P)\left( \tilde{A}-z \right)^{-1}(I-P)\left( A-zI \right)-\left( A-zI \right)P \rbrace=
\]
\[
=\left( \Gamma^{+}\Gamma \right)^{-1}\Gamma^{+}\left\{ 
-A\left( I-P \right)+A(I-P)\left( \tilde{A}-z \right)^{-1}(I-P)\left( A-zI 
\right)-\left( A-zI \right)P \right\}
\]
\[
=\left( \Gamma^{+}\Gamma \right)^{-1}\Gamma^{+}\left\{ 
-\left( A-zI \right)+A(I-P)\left( \tilde{A}-z \right)^{-1}(I-P)\left( A-zI 
\right) \right\}
\]
\[
=\left( \Gamma^{+}\Gamma \right)^{-1}\Gamma^{+}\left\{ 
-I+A(I-P)\left( \tilde{A}-z \right)^{-1}(I-P) \right\}(A-zI).
\]$\square$

Note, if $x_{0}\, \, x_{1}\, \mathellipsis , x_{k-1}$ is a Jordan chain of 
$A$ at the eigenvalue $\alpha \in C$, then it holds 
\[
\, \left( A-zI \right)\left( x_{0}+\left( z-\alpha \right)x_{1}+\mathellipsis 
+\left( z-\alpha \right)^{k-1}x_{k-1} \right)=-\left( z-\alpha \right)^{k}x_{k-1}.
\]
This formula together with (\ref{eq46}) enables us to prove that if $\alpha$ is not a zero of $Q$, then the function $\eta (z):=\hat Q(z)\Gamma^{+}\left( x_{0}+\left( z-\alpha \right)x_{1}+\mathellipsis 
+\left( z-\alpha \right)^{k-1}x_{k-1} \right)= \left( \Gamma^{+}\Gamma \right)^{-1}\Gamma^{+}\left( z-\alpha \right)^{k}x_{k-1}$ is a pole cancellation functions of $Q$ at $\alpha$, cf. \cite[Remark 3.7]{BLu}. 
\\

According to \cite[Proposition 2.1]{Lu1}, for a regular function $Q\in \mathcal{N}_{\kappa 
}(\mathcal{H})$ with representing relation $A$, the inverse $\hat{Q}$ admits 
representation 
\begin{equation}
\label{eq48}
\hat{Q}\left( z \right)=\hat{Q}\left( \bar{z_{0}} \right)+\left( 
z-\bar{z_{0}} \right)\hat{\Gamma }^{+}\left( I+\left( z-z_{0} \right)\left( 
\hat{A}-z \right)^{-1} \right)\hat{\Gamma }
\end{equation}
where $\hat{\Gamma}:=-\Gamma_{z_{0}} {Q(z_{0})}^{-1}$ and it holds 
\begin{equation}
\label{eq50}
\left( \hat{A}-z \right)^{-1}=\left( A-z \right)^{-1}-\Gamma_{z}{Q\left( z 
\right)}^{-1}\Gamma_{\bar{z}}^{+}, \forall z\in \rho \left( A 
\right)\cap \rho \left( \hat{A} \right).
\end{equation}
The following proposition gives us one more relationship between representations 
(\ref{eq36}) and (\ref{eq48}). 

\begin{proposition}\label{proposition12} Let $Q\in \mathcal{N}_{\kappa }(\mathcal{H})$ be holomorphic at $\infty $ and let $Q^{'}(\infty)$ be boundedly invertible. If $\hat{A}$ is the representing linear relation in (\ref{eq48}), then $\hat{A}$ satisfies 
\[
\hat{A}\left( 0 \right)=R\left( P \right)=R(\Gamma).
\]
and $\hat{A}(0)$ is not degenerate. 
\end{proposition}

\textbf{Proof.} Function $Q\in \mathcal{N}_{\kappa }(\mathcal{H})$ that admits representation 
(\ref{eq8}) is a special case of the function that admits representation (\ref{eq2}). 
Let us select a (non-real) point of reference $z_{0}\in \mathcal{D}(Q)\cap \mathcal{D}\left( \hat{Q} 
\right)$, so that $Q\left( z_{0} \right)$ is boundedly invertible. Let us 
introduce $\Gamma_{z_{0}}$ by (\ref{eq28}). Then according to Proposition \ref{proposition2} (ii) 
function $Q$ given by (\ref{eq8}) admits representation (\ref{eq2}) with the same representing 
self-adjoint operator $A$ and ${Q(z_{0})}^{\ast }=\Gamma ^{+}\left( A-\bar{z_{0}} \right)^{-1}\Gamma$. From (\ref{eq50}), for $z=z_{0}$ we get
\begin{equation}
\label{eq52}
\left( \hat{A}-z_{0} \right)^{-1}=\left( A-z_{0} \right)^{-1}-\Gamma 
_{z_{0}}{Q\left( z_{0} \right)}^{-1}\Gamma_{\bar{z_{0}}}^{+}.
\end{equation}
From (\ref{eq28}), it follows
\[
\Gamma_{z_{0}}=\left( A-z_{0} \right)^{-1}\Gamma\,\wedge \Gamma 
_{\bar{z_{0}}}^{+}=\Gamma^{+}\left( A-z_{0} \right)^{-1}.
\]
Substituting this into (\ref{eq52}) gives
\[
\left( \hat{A}-z_{0} \right)^{-1}=\left( A-z_{0} \right)^{-1}-\left( A-z_{0} 
\right)^{-1}\Gamma{Q\left( z_{0} \right)}^{-1}\Gamma^{+}\left( 
A-z_{0} \right)^{-1}
\]
\[
=\left( A-z_{0} \right)^{-1}\left( I-\Gamma{Q\left( z_{0} 
\right)}^{-1}\Gamma^{+}\left( A-z_{0} \right)^{-1} \right).
\]
By substituting here the expression for ${Q\left( z_{0} \right)}^{-1}\Gamma ^{+}$ from (\ref{eq46}) we get
\[
\left( \hat{A}-z_{0} \right)^{-1}=\left( A-z_{0} \right)^{-1}\left( 
I+P\left( -I+A(I-P)\left( \tilde{A}-z_{0} \right)^{-1}(I-P) \right) \right)
\]
\[
=\left( A-z_{0} \right)^{-1}\left( I-P+PA{\left( I-P \right)\left( 
\tilde{A}-z_{0} \right)}^{-1}(I-P) \right).
\]
Hence
\begin{equation}
\label{eq54}
\left( \hat{A}-z_{0} \right)^{-1}=\left( A-z_{0} \right)^{-1}\left( 
I+PA{\left( I-P \right)\left( \tilde{A}-z_{0} \right)}^{-1} \right)\left( 
I-P \right).
\end{equation}
From this we conclude $\ker \left( \hat{A}-z_{0} \right)^{-1}\supseteq R(P)$ 
and, therefore $\hat{A}\left( 0 \right)\supseteq R(\Gamma)$. 

In order to prove $\ker \left( \hat{A}-z_{0} \right)^{-1}\subseteq R\left( 
\Gamma \right)$, assume the contrary, that there exists 

$0\ne (I-P)y\in \ker \left( \hat{A}-z_{0} \right)^{-1}.$ Because, $z_{0}\in 
\rho (A)$ and $A$ is single-valued, from (\ref{eq54}) it follows
\[
\left( I+PA{\left( I-P \right)\left( \tilde{A}-z_{0} \right)}^{-1} 
\right)\left( I-P \right)y=0.
\]
Then, it must be
\[
-\left( I-P \right)y=PA{\left( I-P \right)\left( \tilde{A}-z_{0} 
\right)}^{-1}\left( I-P \right)y=0,
\]
which is a contradiction. Therefore, $\ker \left( \hat{A}-z_{0} 
\right)^{-1}=R(\Gamma)$. \quad $\square$

Note, since the non-real point $z_{0}\in \mathcal{D}(Q)\cap \mathcal{D}(\hat Q)$ was arbitrarily selected, all formulae derived in the proof of Proposition \ref{proposition12} hold for all non-real points $z\in \mathcal{D}(Q) \cap \mathcal{D}(\hat Q)$.

One consequence of Proposition \ref{proposition12} is that function $\hat{Q}$ must have a generalized pole at $\infty$. This means that regular function $\hat{Q}$ does not have a derivative at $\infty$. 

\section{Properties of $\hat{Q}$}\label{s8}

The following theorem is also a consequence of Theorem \ref{theorem4}. 

\begin{theorem}\label{theorem10} Assume that function $Q\in \mathcal{N}_{\kappa}(\mathcal{H})$ is holomorphic at $\infty $, i.e. $ Q \left( z \right):=\Gamma ^{+}\left( A-z 
\right)^{-1} \Gamma $, and assume that operator
\[
Q^{'}\left( \infty\right):=\lim\limits_{z\to \infty}{zQ(z)}
\]
is boundedly invertible. Then for functions 
\begin{equation}
\label{eq62}
\hat{Q}_{1}(z)=\hat{S}+z\hat{G} \in \mathcal{N}_{\kappa_{1}}\left( \mathcal{H} \right),
\end{equation}
and
\begin{equation}
\label{eq63}
\hat{Q}_{2}\left( z \right):=\tilde{\Gamma }^{+}\left( \tilde{A}-z 
\right)^{-1}\tilde{\Gamma }\in \mathcal{N}_{\kappa_{2}}\left( \mathcal{H} \right), 
\end{equation}
where operators $\hat S$, $\hat G$ and $\tilde \Gamma $ are given by equations (\ref{eq42a}) and (\ref{eq42b}), the inverse function $\hat{Q}\left( z \right)$ has decomposition 
\begin{equation}
\label{eq64}
\hat{Q}\left( z \right)=\hat{Q}_{1}\left( z \right)+\hat{Q}_{2}\left( z \right).
\end{equation}
That decomposition has the following properties:

\begin{enumerate}[(i)]
\item It must be $\hat{Q}_{1}\not\equiv 0$\, while function $\hat{Q}_{2}$ may be zero function in some cases. $\hat{Q}_{1}$ has only one generalized pole, it is at $\infty $, while $\hat{Q}_{2}$ is holomorphic at $\infty $.
\item Finite generalized zeros of $Q$, coincide with generalized poles of $\hat{Q}_{2}$ including multiplicities. 
\item $\hat{Q}_{1}\in \mathcal{N}_{\kappa_{1}}(\mathcal{H})$, where negative index $\kappa_{1}$ is equal to the number of negative eigenvalues of the bounded self-adjoint operator $-Q^{'}\left( \infty\right)$ in the Hilbert space $\mathcal{H}$ and that is equal to negative index of $P\mathcal{K}$.
\item $\kappa_{1}+\kappa_{2}=\kappa $.
\end{enumerate}
\end{theorem}

\textbf{Proof.} (i) According to above definitions of $\hat Q_{1}$ and $\hat Q_{2}$, and  (\ref{eq42}), it holds $\hat{Q}\left( z \right)=\hat{Q}_{1}\left( z \right)+\hat{Q}_{2}\left( z \right)$. According to Proposition \ref{proposition12}, $\hat{Q}$ has generalized pole at $\infty $. Since 
representing operator $\tilde{A}$ of $\hat{Q}_{2}$ is bounded operator, according to Lemma \ref{lemma3}  $\hat{Q}_{2}$ is holomorphic at $\infty $. Therefore, $\hat{Q}_{1}\not\equiv 0$ and it must have generalized pole at $\infty $. 
According to Example \ref{example2} it is possible to have $\hat{Q}_{2}\equiv 0$.

(ii) The statement follows immediately from (i) and formula (\ref{eq64}). 

(iii) Note, representation (\ref{eq62}) of $\hat{Q}_{1}$ is not a typical operator 
representation of a generalized Nevanlinna function, because $A-zI$ is not a 
resolvent. 

We know $\hat{Q}\in \mathcal{N}_{\kappa }(\mathcal{H})$ and $\kappa_{1}+\kappa_{2}\ge \kappa 
$. Let us denote by $\kappa^{'}$ and $\kappa^{''}$ negative indexes of 
subspaces $P\mathcal{K}$ and $\left( I-P \right)\mathcal{K}$, respectively. Then, according to 
(\ref{eq32}) $\kappa^{'}+\kappa^{''}=\kappa $. 

For any $f,\, g\in \mathcal{H}$ we have
\[
\left( \frac{\hat{Q}_{1}\left( z \right)-{\hat{Q}_{1}\left( w \right)}^{\ast 
}}{z-\bar{w}}f,g \right)=\left( \left( \Gamma^{+}\Gamma 
\right)^{-1}f,g \right).
\]
Hence, $\kappa_{1}$ equals number of negative eigenvalues of $\left( \Gamma^{+}\Gamma \right)^{-1}$. Since $\left( \Gamma^{+}\Gamma 
\right)^{-1}$ is bounded, hence defined on the whole $\mathcal{H}$, we can consider 
$f=\Gamma^{+} \Gamma f_{0}$ and $g=\Gamma^{+} \Gamma g_{0} $, 
where $f_{0}$ and $g_{0}$ run through entire $\mathcal{H}$ when $f$ and $g$ run through $\mathcal{H}$. 
Therefore 
\[
\left( \left( \Gamma^{+} \Gamma \right)^{-1} f, g \right)=\left[ \Gamma f_{0} , \Gamma g_{0} \right].
\]
Because $R\left( \Gamma \right)=R(P)$, we conclude that $\kappa_{1}=\kappa 
^{'}$. Real number $\alpha <0$ is an eigenvalue of $\Gamma^{+}\Gamma=-Q^{'}\left( \infty\right)$ if and only if $\alpha^{-1}<0$ is an eigenvalue of $\left( \Gamma^{+}\Gamma \right)^{-1}$. Hence, statement (iii) follows.

(iv)
\[
\kappa_{1}=\kappa^{'}\Rightarrow \kappa^{'}+\kappa_{2}\ge \kappa =\kappa 
^{'}+\kappa^{''}\Rightarrow \kappa_{2}\ge \kappa^{''}
\]
Because $\tilde{A}$, the representing operator of $\hat{Q}_{2}$, is self-adjoint 
operator in $\left( I-P \right)\mathcal{K}$, it must be $\kappa_{2}\le \kappa^{''}$. 
Therefore, $\kappa_{2}=\kappa^{''}$ and
\[
\kappa_{1}+\kappa_{2}=\kappa . 
\]
That proves iv).\quad $\square$
\\

In the following example we will show how Theorem \ref{theorem10} can be applied to a 
concrete generalized Nevanlinna functions. 

\begin{example}\label{example4}Let 
\[
Q\left( z \right)=\left[ {\begin{array}{*{20}c}
\frac{-(1+z)}{z^{2}} & \frac{1}{z}\\[1mm]
\frac{1}{z} & \frac{1}{1+z}\\[1mm]
\end{array} } \right].
\]
\end{example}
The function $Q$ has representation (\ref{eq8})
\[
Q\left( z \right)=\Gamma^{+}\left( A\mathrm{-}z \right)^{\mathrm{-1}}\Gamma ,
\]
where the space $K=C^{3}$. In that representation fundamental symmetry, and representing operators of $Q$ are:
\[
J=\left[ {\begin{array}{*{20}c}
0 & 1 & 0\\
1 & 0 & 0\\
0 & 0 & -1\\
\end{array} } \right],\, \, A=\left[ {\begin{array}{*{20}c}
0 & 1 & 0\\
0 & 0 & 0\\
0 & 0 & -1\\
\end{array} } \right],\, \, \Gamma =\left[ {\begin{array}{*{20}c}
0.5 & -1\\
1 & 0\\
0 & -1\\
\end{array} } \right],\, \, \Gamma^{+}=\Gamma^{\ast }J=\left[ 
{\begin{array}{*{20}c}
1 & 0.5 & 0\\
0 & -1 & 1\\
\end{array} } \right].
\]
Here, $\Gamma^{\ast }:C^{3}\to C^{2}$ is adjoint operator of $\Gamma $ with 
respect to Hilbert spaces $C^{2}$ and $C^{3}$. It is easy to see that this 
representation is minimal. From the shape of the fundamental symmetry $J$ we 
conclude $\kappa =2$, i.e. $Q\in \mathcal{N}_{2}(C^{2})$. We have 
\[
\hat{Q}\left( z \right)=\left[ {\begin{array}{*{20}c}
\frac{z^{2}}{2(1+z)} & -\frac{z}{2}\\
-\frac{z}{2} & \frac{-(1+z)}{2}\\
\end{array} } \right]\in N_{2}(C^{2}).
\]
Limit (\ref{eq30b}) gives
\[
\Gamma^{+}\Gamma =\left[ {\begin{array}{*{20}c}
1 & -1\\
-1 & -1\\
\end{array} } \right],\, \left( \Gamma^{+}\Gamma \right)^{-1}=\left[ 
{\begin{array}{*{20}c}
0.5 & -0.5\\
-0.5 & -0.5\\
\end{array} } \right].
\]
This means that conditions of Theorem \ref{theorem10} are satisfied. 

Let us calculate $\hat{Q}_{1}\left( z \right)$. By substituting matrices 
$\left( \Gamma^{+}\Gamma \right)^{-1},\, \, \Gamma^{+},\, \, \Gamma \, $ 
into formulae for $\hat G$ and $\hat S$, we obtain   
\[
\hat{Q}_{1}\left( z \right)=\left[ {\begin{array}{*{20}c}
\frac{-1+z}{2} & -\frac{z}{2}\\
-\frac{z}{2} & -\frac{1+z}{2}\\
\end{array} } \right].
\]

Let us now find $\hat{Q}_{2}\left( z \right)$ by means of formulae 
(\ref{eq63}). In order to do that, we have first to find matrices for projections 
$P$ and $(I-P)$. By means of formula (\ref{eq31}) we get 
\[
P=\left[ {\begin{array}{*{20}c}
0.75 & 0.125 & 0.25\\
0.5 & 0.75 & -0.5\\
0.5 & -0.25 & 0.5\\
\end{array} } \right],\, \, I-P=\, \left[ {\begin{array}{*{20}c}
0.25 & -0.125 & -0.25\\
-0.5 & 0.25 & 0.5\\
-0.5 & 0.25 & 0.5\\
\end{array} } \right].
\]
Obviously, range $\left( I-P \right)=1$, i.e. $\dim {\left( I-P \right)K=1}$. We also have 
\[
\left( I-P \right)A\left( I-P \right)-z\left( I-P \right)=\, \left[ 
{\begin{array}{*{20}c}
-0.25 & 0.125 & 0.25\\
0.5 & -0.25 & -0.5\\
0.5 & -0.25 & -0.5\\
\end{array} } \right]-z\left[ {\begin{array}{*{20}c}
0.25 & -0.125 & -0.25\\
-0.5 & 0.25 & 0.5\\
-0.5 & 0.25 & 0.5\\
\end{array} } \right],
\]
\[
\tilde{\Gamma }:=\left( I-P \right)A\Gamma \left( \Gamma^{+}\Gamma 
\right)^{-1}=\left[ {\begin{array}{*{20}c}
0.25 & 0\\
-0.5 & 0\\
-0.5 & 0\\
\end{array} } \right],\, \, \tilde{\Gamma }^{+}=\tilde{\Gamma 
}^{\ast }J=\left[ {\begin{array}{*{20}c}
-0.5 & 0.25 & 0.5\\
0 & 0 & 0\\
\end{array} } \right].
\]
Obviously, $\tilde{\Gamma }$, and $\tilde{\Gamma }^{+}$, each have only one linearly independent row,  column, respectively. Therefore, 
\textbf{operators} $\tilde{\Gamma },\, \tilde{\Gamma }^{+}$ can be 
represented by equivalent matrices, i.e. we can write
\[
\tilde{\Gamma }:=\left[ {\begin{array}{*{20}c}
0.25 & 0\\
0 & 0\\
0 & 0\\
\end{array} } \right],\, \, \tilde{\Gamma }^{+}=\left[ 
{\begin{array}{*{20}c}
-0.5 & 0 & 0\\
0 & 0 & 0\\
\end{array} } \right].
\]
Accordingly, we will write in the equivalent matrix form
\[
\left( I-P \right)A\left( I-P \right)-z\left( I-P \right)=\, \left[ 
{\begin{array}{*{20}c}
-0.25-0.25z & 0 & 0\\
0 & 0 & 0\\
0 & 0 & 0\\
\end{array} } \right].
\]
Then, the matrix form of the \textbf{operator }
\[
\left( I-P \right)\left( \tilde{A}-z \right)^{-1}\left( I-P \right)=\left( 
{\begin{array}{*{20}c}
\left( \tilde{A}-z \right)^{-1} & 0\\
0 & 0\\
\end{array} } \right)
\]
is
\[
\left[ {\begin{array}{*{20}c}
\frac{-4}{1+z} & 0 & 0\\
0 & 0 & 0\\
0 & 0 & 0\\
\end{array} } \right].
\]
Now, according to (\ref{eq63}) we calculate 
\[
\hat{Q}_{2}\left( z \right):=\tilde{\Gamma }^{+}\left( \tilde{A}-z 
\right)^{-1}\tilde{\Gamma }=\left[ {\begin{array}{*{20}c}
-0.5 & 0 & 0\\
0 & 0 & 0\\
\end{array} } \right]\left[ {\begin{array}{*{20}c}
\frac{-4}{1+z} & 0 & 0\\
0 & 0 & 0\\
0 & 0 & 0\\
\end{array} } \right]\left[ {\begin{array}{*{20}c}
0.25 & 0\\
0 & 0\\
0 & 0\\
\end{array} } \right].
\]
Thus
\[
\hat{Q}_{2}\left( z \right)=\left[ {\begin{array}{*{20}c}
\frac{1}{2(1+z)} & 0\\
0 & 0\\
\end{array} } \right].
\]
We obtained the decomposition (\ref{eq64}) of $\hat{Q}\left( z \right)$:
\[
\left[ {\begin{array}{*{20}c}
\frac{z^{2}}{2(1+z)} & -\frac{z}{2}\\
-\frac{z}{2} & \frac{-(1+z)}{2}\\
\end{array} } \right]=\left[ {\begin{array}{*{20}c}
\frac{-1+z}{2} & -\frac{z}{2}\\
-\frac{z}{2} & -\frac{1+z}{2}\\
\end{array} } \right]+\left[ {\begin{array}{*{20}c}
\frac{1}{2(1+z)} & 0\\
0 & 0\\
\end{array} } \right].
\]

There are many decompositions of the function $\hat{Q}$. For this decomposition, we know that the  following claims hold: 

- Because Hermitian matrix $\Gamma^{+}\Gamma $ has one simple negative 
eigenvalue, according to Theorem \ref{theorem10} (iii) the function 
$\hat{Q}_{1}$ has negative index $\kappa_{1}=1.$ 

- Because, $\kappa =2$, according to Theorem \ref{theorem10} (iv), it must be $\kappa_{2}=1$. 

- According to Theorem \ref{theorem10} (ii), $z=-1$ is zero of the function $Q$. Indeed, it is a pole of $\hat Q_{2}$ with pole cancellation function $\eta (z) =\left[ {\begin{array}{*{20}c}
1+z\\
0\\
\end{array} } \right] $, according to \cite[Definition 3.1]{BLu}. $\square$
\\

In this example we have demonstrated how to use formulae given in Theorem \ref{theorem10} to obtain decomposition (\ref{eq64}). The example was 
selected to be as simple as possible to make it readable. In 
more complicated cases, the calculation of 
\[
\hat{Q}_{1}(z)=\hat{S}+z\hat{G}
\]
remains  simple, while calculation of $\hat{Q}_{2}\left( z \right)$ can get very involved . 

Fortunately, Theorem \ref{theorem10} enables us to avoid the difficult calculation of $\hat{Q}_{2}$ given by formula (\ref{eq63}). Instead, we can obtain $\hat{Q}_{2}$ by formula $\hat{Q}_{2}\left( z \right):=\hat{Q}\left( z \right)-\hat{Q}_{1}\left( z \right)$.

In general case, it is an interesting task to decompose a generalized Nevanlinna function into a sum that preserves the number of negative squares, i.e. $Q =Q_{1}+Q_{2}$ and $\kappa =\kappa_{1}+\kappa_{2}$.

Street Address: 

BML, Actuarial department

120 Royall St. Canton MA, 02021 

Emails: 

Muhamed.Borogovac@gmail.com

Muhamed$_{-}$Borogovac@bostonmutual.com 
\maketitle
\tableofcontents

\parskip=4pt
\parindent=0cm
\end{document}